\documentclass[12pt,reqno]{amsart} 
\usepackage[all,2cell,arrow]{xy}
\usepackage[T1]{fontenc}
\usepackage{amsmath}
\usepackage{amsthm}
\usepackage{amsfonts}
\usepackage{amssymb}
\usepackage[all]{xy}
\usepackage{hyperref}
 \usepackage{color}

\newcommand{\ba}{\begin{array}}
\newcommand{\ben}{\begin{enumerate}}
\newcommand{\bt}{\begin{tabular}}

\newcommand{\ea}{\end{array}}
\newcommand{\ed}{\end{diagram}}
\newcommand{\een}{\end{enumerate}}
\newcommand{\et}{\end{tabular}}

\newcommand{\bit}{\begin{itemize}}
\newcommand{\eit}{\end{itemize}}

\newenvironment{theo.}{\removelastskip\par\medskip
\refstepcounter{subsection}
\noindent{\bf\thesubsection. Th�or�me.}\em}{\par\em\medskip}

\newenvironment{prop.}{\removelastskip\par\medskip
\refstepcounter{subsection}
\noindent{\bf\thesubsection. Proposition.}\em}{\par\em\medskip}

\newenvironment{coro.}{\removelastskip\par\medskip
\refstepcounter{subsection}
\noindent{\bf\thesubsection. Corollaire.}\em}{\par\em\medskip}

\newenvironment{defi.}{\removelastskip\par\medskip
\refstepcounter{subsection}
\noindent{\bf \thesubsection. D�finition.}}{\par\medskip}

\newcounter{subsec}[section]

\newcommand{\monr}{\ar@{{^ (}->}[r]}
\newcommand{\monl}{\ar@{{^ (}->}[l]}
\newcommand{\monu}{\ar@{{^ (}->}[u]}
\newcommand{\mond}{\ar@{{^ (}->}[d]}
\newcommand{\monrr}{\ar@{{^ (}->}[rr]}
\newcommand{\monll}{\ar@{{^ (}->}[ll]}
\newcommand{\monuu}{\ar@{{^ (}->}[uu]}
\newcommand{\mondd}{\ar@{{^ (}->}[dd]}
\newcommand{\mondr}{\ar@{{^ (}->}[dr]}
\newcommand{\monur}{\ar@{{^ (}->}[ur]}
\renewcommand{\to}{\longrightarrow}

\renewcommand{\ker}{\ensuremath{\mathsf{Ker\,}}}

\newcommand{\Ab}{\mathsf{Ab}}
\newcommand{\Gp}{\mathsf{Gp}}

\newcommand{\Ext}{\ensuremath{\mathsf{Ext(\mathbb A})}}

\newcommand{\CExtF}{\ensuremath{\mathsf{{CExt}_{\mathbb B}(\mathbb A)}}}

\newcommand{\Ac}{\ensuremath{\mathbb{A}}}

\newcommand{\Bc}{\ensuremath{\mathbb{B}}}

\newbox\pullbackbox
\setbox\pullbackbox=\hbox{\xy \POS(4,0)\ar@{-} (0,0) \ar@{-} (4,4)
\endxy}

\newbox\pushoutbox
\setbox\pushoutbox=\hbox{\xy \POS(0,4)\ar@{-} (0,0) \ar@{-} (4,4)
\endxy}

\begin{document}

\title{Some aspects of semi-abelian homology and protoadditive functors}

\author{ \small{\emph{Tomas Everaert and Marino Gran}  }}
\address{Institut de Recherche en Math\'ematique et  Physique,
Universit\'e catholique de Louvain, Chemin du Cyclotron 2, Louvain-la-Neuve,
Belgium}

\begin{abstract}{In this note some recent developments in the study of homology in semi-abelian categories are briefly presented. In particular the role of protoadditive functors in the study of Hopf formulae for homology is explained. }
\end{abstract}
\vspace{2cm}

\maketitle

The discovery of higher Hopf formulae for the homology of a group, due to Ronald Brown and Graham Ellis \cite{BE}, has naturally led to some new perspectives in non-abelian homological algebra. An important advance in this area was made by George Janelidze, who found a connection between group homology and categorical Galois theory \cite{Jan, J1, J2, BoJ}, the latter being a wide extension of Alexander Grothendieck's theory \cite{Gro}. Among other things, this cleared the path for the discovery of higher Hopf formulae for the homology of general algebraic structures \cite{EVL1, EVL, EvG, EGV, Ev, EG, DEG, Proto}. Here, a crucial role is played by the so-called \emph{higher order central extensions}, which are the covering morphisms with respect to certain Galois structures induced by a reflection 
\begin{equation}\label{adjunction}
\xymatrix{ {\mathbb B \, } \ar@<-1ex>[r]_-{U} & {\, \mathbb A, \, }
\ar@<-1ex>[l]^-{^{\perp}}_-{I}  }
\end{equation}
whose left adjoint $I \colon \mathbb A \rightarrow \mathbb B$ is sometimes called the ``coefficient functor''. Here, $\mathbb A$ could, for instance, be the variety of groups, $\mathbb B$ its subvariety of abelian groups, and $I$ the abelianisation functor. In this case, the induced higher order central extensions are related to the Brown-Ellis Hopf-formulae, as explained below. More generally, higher order central extensions can be defined for any semi-abelian category $\mathbb A$ (\cite{JMT} e.g. the varieties of groups, rings, Lie algebras, (pre)crossed modules, compact groups, or any abelian category) and any Birkhoff subcategory  (i.e. a reflective subcategory closed under subobjects and regular quotients) $\mathbb B$ of $\mathbb A$. When, moreover,  $\mathbb A$ has enough projectives, one obtains higher Hopf formulae for the homology induced by the reflector (or, coefficient functor) $I\colon \mathbb A\to \mathbb B$. 

In order to explicitly determine the Hopf-Brown-Ellis formulae for homology in specific algebraic contexts,  it is crucial to find suitable descriptions of the higher central extensions, as for instance in terms of algebraic conditions using ``generalised commutators''. In general, this is a non-trivial problem, about which we are going to say more in what follows. 

Note that, in fact, different approaches to obtaining higher Hopf formulae exist, which can be used in different categorical contexts, based on the comonadic homology theory of Barr and Beck \cite{BB, EGV, EVL1, EVL}, on the abstract notion of Galois group \cite{JanHopf, DEG, Mathieu}, or on the theory of satellites \cite{GV}. These methods essentially coincide in the situation described above, namely for $\mathbb A$ a semi-abelian category with enough projectives and $\mathbb B$ a Birkhoff subcategory of $\mathbb A$. 
      
Assume that we are in this situation. In this case, the reflector $I\colon \mathbb A\to \mathbb B$ induces a first ``centralisation functor'' $I_1$ from the category $\Ext $ of extensions (i.e regular epimorphisms) in $\Ac$  to the full subcategory $\CExtF$ of extensions that are central with respect to $\Bc$:
 \begin{equation}\label{adjcentr}
\xymatrix@=30pt{\CExtF   \ar@<-1ex>[r]_-{U_1} & {\Ext .} \ar@<-1ex>[l]^-{^{\perp}}_-{I_1} }
\end{equation}
This functor $I_1$ is the left adjoint of the inclusion functor $U_1$. The notion of centrality comes from categorical Galois theory and depends on the choice of the Birkhoff subcategory $\Bc$. It is defined in purely categorical terms \cite{JK, EGV}.

The centralisation $I_1 (f)$ of an extension $f$ is given by a quotient
$$
\xymatrix{A \ar[rr]^-{{\eta}_f} \ar[dr]_{f}& & A/[\ker(f), A]_{\Bc} \ar[dl]^{I_1 (f)} \\& B, &}
$$
where $[\ker(f), A]_{\Bc}$ may be thought of as a commutator of $\ker (f)$ and $A$, defined relatively to $\mathbb B$. For instance, in the classical case of the reflection
$$
\xymatrix@=30pt{
{\Ab \, } \ar@<-1ex>[r]_-{U} & {\, \mathsf{Gp} \, }
\ar@<-1ex>[l]^-{^{\perp}}_-{\mathsf{ab}} }
$$
this relative commutator is simply the group-theoretical commutator of normal subgroups: $$[\ker(f),A]_{\Ab}= [\ker(f), A].$$  Hence, in this case $I_1$ is the usual centralisation functor from the category of group extensions to its full subcategory of central extensions in the classical sense.

Remark that the commutator $[\ker(f), A]$ appears as the denominator in the Hopf formula for the second integral homology group: given a free presentation 
$$\xymatrix{0 \ar[r]& K \ar[r] & F \ar[r] & G \ar[r] & 0,}
$$
of a group $G$, Hopf's formula tells us that the second homology group is given by the quotient
$$
H_2 (G, {\mathbb Z})= \frac{K \cap [F, F]}{[K, F]}.
$$
This is not a coincidence: a similar phenomenon occurs for the higher-order homology groups, where the subgroups appearing as denominator of the Brown-Ellis Hopf formulae are exactly what is required to transform higher-order extensions into a higher-order \emph{central} extensions, universally.

To illustrate this idea, consider the case of the \emph{third} homology group. For this, let 
\begin{equation}\label{doppia}
\xymatrix{ 
F \ar[r] \ar[d] & F/K_2 \ar@{.>}@{->}[d] \\
 F/K_1 \ar@{->}[r] & G\cong F/K_1 \cdot K_2}
\end{equation}
 be a double presentation of a group
$G$ so that $F$, $F/K_1$ and $F/K_2$ are free groups and the square is a \emph{double extension}: a pushout of surjective homomorphisms. 
As shown by Brown and Ellis, the third integral homology group $H_3 (G, {\mathbb Z})$ ($=H_3(G,{ \mathsf{Ab}})$)
of $G$ is given by
\begin{equation}\label{formule3}
H_3 (G, {\mathbb Z}) = \frac{[F,F]\cap K_1 \cap K_2}{[K_1 \cap K_2,F]\cdot  [K_1, K_2]}.
\end{equation}
Once more, the denominator gives precisely the normal subgroup by which one has to ``quotient out'' the group $F$ in order to make this double extension a  double \emph{central} extension, universally:
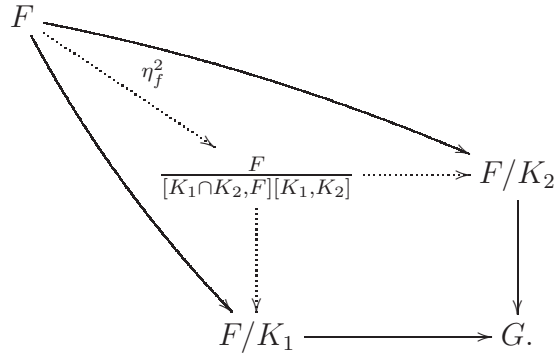
\[
\vcenter{\xymatrix@=40pt{ 
F  \ar@/^/@{->}[drr]^{} \ar@/_/@{->}[drd]_{} \ar@{.>}[rd]^{\eta^2_f} \\
& \frac{F}{[K_1 \cap K_2, F][K_1, K_2]}  \ar@{.>}[r] \ar@{.>}[d] & F/K_{2} \ar@{.>}@{->}[d]^{} \\
& F/K_1 \ar@{->}[r]_{} & G.}}
\]
Once again, the notion of centrality comes from categorical Galois theory, and this time depends on the induced reflection \eqref{adjcentr}. 

Now, the formula \eqref{formule3} is a special case of the general Hopf formula for the third homology corresponding to a reflection \eqref{adjunction}, with $\mathbb A$ an arbitrary semi-abelian category with enough projectives, and $\mathbb B$ any Birkhoff subcategory of $\mathbb A$ \cite{EGV}: starting from a double presentation $f$ of the form \eqref{doppia}, with $F$, $F/K_1$ and $F/K_2$ regular projective objects of $\Ac$, the third homology object is given by a quotient
 \begin{equation}\label{hopfg}
H_3 (G, \Bc) = \frac {[F,F]_{\Bc} \cap K_1 \cap K_2}{L_2[f]_{\Bc}}.
\end{equation}
As in the case of groups, {also in a general semi-abelian category $\mathbb A$} the denominator $L_2[f]_{\Bc}$ of this generalised Hopf formula {relative to $\mathbb B$} is the normal subobject of $F$ that has to be ``quotiented out'' in order to universally turn the double extension \eqref{doppia} into a double central extension. Hence, in particular, for $\mathbb A=\Gp$ and $\mathbb B=\Ab$, we have the equality
$$L_2[f]_{\Ab}= [K_1 \cap K_2,F] \cdot [K_1, K_2],$$
and the formula \eqref{formule3} appears as a special case of \eqref{hopfg}. In general, for a given Birkhoff reflection \eqref{adjunction} in a semi-abelian category, $L_2[f]_{\Bc}$ may be difficult to compute.

Similar formulae exist for the higher homology objects, again valid in any semi-abelian category $\Ac$ with enough projectives and for any Birkhoff subcategory $\Bc$ of $\Ac$. This yields, at least in principle, a description of all homology objects $H_n (A, \Bc)$ ($n\geq 2$). In practice, a suitable characterisation of the higher-order central extensions is required.

In some cases it has been possible to compute these formulae explicitly. For example, this has been done in \cite{EGV} for $\Ac$ the variety of precrossed modules and $\Bc$ its subvariety of crossed modules, or for $\Ac$ the variety of groups and $\Bc$ its subvariety of nilpotent groups of a fixed class $k\ge 1$ (see \cite{Don}), or the variety of solvable groups of a fixed class $k\ge 1$. Similar results have been obtained in the categories of Leibniz and of Lie $n$-algebras in \cite{CKLV}.

However, in general, computing the Hopf formulae explicitly is a non-trivial task. One possible strategy is to look for suitable conditions on the coefficient functor $I \colon \Ac \rightarrow \Bc$ that facilitate such computations. In \cite{EG, Proto} we have shown that a natural such condition is the requirement that the reflector $I$ is a \emph{protoadditive functor}. This notion extends the one of additive functor to the non-additive context of pointed protomodular categories \cite{Bourn}: when $\Ac$ and $\Bc$ are pointed protomodular categories (for instance, $\mathbb A$ and $\mathbb B$ could be semi-abelian), a functor $I \colon \Ac \rightarrow \Bc$ is \emph{protoadditive} if, for any split short exact sequence 
\[
\xymatrix{0 \ar[r]& K \ar[r]^k & A \ar@<-1 ex>[r]_f & B \ar@<-1ex>[l]_s \ar[r] &0 }
\]
in $\Ac$ (i.e. $f\circ s=1_B$ and $k=\ker (f)$), its image 
$$\xymatrix{0 \ar[r]& I(K) \ar[r]^{I(k)} & I(A) \ar@<-1 ex>[r]_{I(f)} & I(B)  \ar@<-1ex>[l]_{I(s)} \ar[r] &0 }$$
under $I$ is a split short exact sequence in $\Bc$.
Whenever the coefficient functor $I \colon \mathbb A \rightarrow \mathbb B$ is protoadditive, explicit Hopf formulae can be established in different algebraic and topological contexts. In particular, the protoadditivity condition is fundamental to explore some new Galois theories induced by torsion theories (see also \cite{BG, GR, GJ}).

We refer the interested reader to the articles \cite{EG} and \cite{Proto} for a thorough study of the theory of protoadditive functors and their use in semi-abelian homological algebra. In the first article we study the homology of $n$-fold internal groupoids in a semi-abelian category: these results apply in particular to the so-called $\mathsf{cat}^n$-groups in the sense of \cite{Lo}. The crucial point there is that the connected components functor $\pi_0 \colon \mathsf{Grpd}(\Ac) \rightarrow \Ac$ is protoadditive whenever $\Ac$ is semi-abelian. In \cite{Proto} the general theory of protoadditive functors is investigated, as well as the so-called \emph{derived torsion theories} of a torsion theory having the reflector to its torsion-free subcategory \emph{protoadditive}. These derived torsion theories induce a chain of Galois structures in the categories of higher order extensions. The results concerning the homology objects can be applied in particular to some new torsion theories in the category of compact groups and, more generally, in any category of compact semi-abelian algebras introduced in \cite{BoC}. Further developments in this direction, including some new results in the categories of commutative rings and of topological groups, for instance, can be found in \cite{DEG, Mathieu}.
\vspace{2mm}

{{\bf{Acknowledgement}.}
This article is partly based on the text of a talk given by the second author on this joint work at the \emph{S\'eminaire Itin\'erant des Cat\'egories} that took place} at the Universit\'e Paris Did\'erot on October 25, 2009.

\end{document}